\documentclass[11pt]{article}
\usepackage{amscd, amsmath, amssymb}

\title{The automorphism group of an affine quadric}
\author{Burt Totaro}
\date{  }

\def\R{\text{\bf R}}

\def\P{\text{\bf P}}

\def\arrow{\rightarrow}

\def\qed{\ QED \vspace{.1in}}

\def\Aut{\text{Aut}}

\begin{document}
\maketitle

\newtheorem{theorem}{Theorem}[section]
\newtheorem{corollary}[theorem]{Corollary}
\newtheorem{lemma}[theorem]{Lemma}
\newtheorem{conjecture}[theorem]{Conjecture}

We determine the automorphism group for
a large class of affine quadrics over a field, viewed
as affine algebraic varieties. The proof uses a fundamental theorem
of Karpenko's in the theory of quadratic forms \cite{Karpenko},
along with some useful arguments of birational geometry.
In particular, we find that the automorphism
group of the $n$-sphere $\{x_0^2+\cdots + x_n^2=1\}$ over the real numbers
is just the orthogonal
group $O(n+1)$ whenever $n$ is a power of 2. It is not known
whether the same is true for arbitrary $n$.
This result is reminiscent of Wood's theorem that
when $n$ is a power of 2, every real polynomial mapping from the
$n$-sphere to a lower-dimensional sphere is constant \cite{Wood}.

The background for these results is that almost all geometric
tools work better for projective varieties than for affine varieties,
because of the lack of compactness.
Even basic questions like determining the
automorphism group of an affine variety, or whether two
affine varieties are isomorphic, can be difficult.
Of course, some cases are easy. Consider an affine variety
$X-D$ where $X$ is projective. If
$X$ is of general type, or more generally if the pair $(X,D)$ is
of log-general type in Iitaka's sense (for example,
when $D$ is a smooth hypersurface of
degree at least $n+2$ in $X=\P^n$), then the affine variety $X-D$
has finite automorphism group \cite[Theorem 11.12]{Iitaka}.
But when both $X$ and $D$ are of low degree in some sense,
then the automorphism group of $X-D$ is not at all
understood. This justifies studying the basic case of affine
quadrics.

The key ingredient of the proof of the main Theorem \ref{witt}
is that an anisotropic projective quadric with
first Witt index equal to 1 is not ruled. More generally,
a fundamental problem of birational geometry
is to determine which varieties over a field are ruled.
For example, Koll\'ar proved that a large class
of rationally connected complex hypersurfaces are non-rational
by showing that they are not even ruled
\cite[Theorem V.5.14]{Kollarbook}.
For anisotropic quadrics over a field, we give a conjectural
answer to the problem of ruledness: they should be 
ruled if and only if the first Witt index is greater than 1
(Conjecture \ref{quadruled}).
Section \ref{questions} gives some evidence: in particular, the conjecture
is true for quadratic forms of dimension at most 9 (thus for projective
quadrics of dimension at most 7).

\section{Infinite-dimensional automorphism groups}

In this section, we show that all affine quadrics of dimension
at least 2 whose homogeneous
part of degree 2 is isotropic
have infinite-dimensional automorphism group. Here an affine quadric
is the affine scheme defined by the vanishing of a polynomial $f$
of degree at most 2 over a field $k$.
We assume only that $f$ is not a constant;
we do not assume any nondegeneracy, and the field $k$ may be
of any characteristic including 2.
A quadratic form is defined to be
a homogeneous polynomial
of degree 2 on a finite-dimensional vector space $V$
over a field $k$. We say that a quadratic form $q$ is isotropic
if there is a nonzero vector $x$ in $V$ with $q(x)= 0$.

\begin{lemma}
\label{inf}
Let $X=\{f=0\}$ be an affine quadric in an $n$-dimensional vector space
over a field $k$
with $n\geq 3$. If the homogeneous part of degree 2 of $f$ is isotropic,
then the automorphism group of $X$ is infinite-dimensional.
Also, for any isotropic quadratic form $q$
of dimension $n\geq 3$ over a field $k$,
the complement $\P^{n-1}-\{q = 0\}$ has infinite-dimensional automorphism
group.
\end{lemma}

When $n=3$, more is known. For 
an isotropic affine quadric surface $xy+az^2=b$
over a field, or the complement of an isotropic conic $xy+az^2=0$
in $\P^2$, Gizatullin and Danilov computed the automorphism
groups explicitly, as amalgamated free products \cite{GD2}.

For the purposes of this section, we will say that the automorphism
group of an algebraic variety $X$ over a field $k$ is infinite-dimensional
if there are unirational varieties $S$ over $k$ of arbitrarily
large dimension which have a  morphism $S\times X\arrow X$ that
induces an injection of $S(k)$ into $\Aut(X)$. (Here unirationality,
being the image of a dominant rational map from projective space,
ensures that the set $S(k)$ of $k$-rational points
is ``big'', at least for infinite
fields $k$.) One could consider weaker definitions, but Lemma
\ref{inf} will show that certain affine quadrics have infinite-dimensional
automorphism group over $k$ in this strong sense.

{\bf Proof. }
We can view $X$ as a projective quadric minus its intersection
with a hyperplane. Let the projective quadric
be defined by a quadratic form $q$ on
a vector space $V$, while the hyperplane corresponds to a hyperplane
$W\subset V$. We are assuming that $q$ restricted to $W$ is isotropic,
and so we can choose a nonzero vector $x$ in $W$ with $q(x)=0$.
The main point is to exhibit a nontrivial homomorphism
from the additive group $G_a$ to the orthogonal group $O(W\subset V)$
of linear automorphisms of $V$ that preserve $q$ and the subspace $W$.

Let $\langle u,v\rangle =q(u+v)-q(u)-q(v)$ be the symmetric
bilinear form associated to $q$. If the isotropic vector $x$ in $W$
is orthogonal to all of $V$,
we choose a nonzero linear form $g:V\arrow k$ that vanishes on $x$.
Then 
$$\varphi_t(z)=z+tg(z)x$$
is a nontrivial homomorphism from $G_a$ to $O(W\subset V)$. Next,
suppose $x$ is not orthogonal to all of $V$. Let $y$ be a vector
in $W$ orthogonal to $x$ and not a multiple of $x$; this exists,
since $W$ has dimension at least 3. Then the Siegel transvection
\cite[III.1.5]{Baeza}
$$\varphi_t(z)=z+t\langle z,x\rangle y-t\langle z,y\rangle x
-t^2q(y)\langle z,x\rangle x$$
is a nontrivial homomorphism from $G_a$ to $O(W\subset V)$.

Thus, we have a nontrivial action of $G_a$ on the affine quadric $X$. This 
gives a injective homomorphism from the $k$-vector space of $G_a$-invariant
regular functions $f$ on $X$ to the automorphism group of $X$, given by
$z\mapsto \varphi_{f(z)}(z)$. Since
$X$ has dimension at least 2, we see by hand that the vector space
of $G_a$-invariant functions on $X$ is infinite-dimensional
(use arbitrary polynomials in the coordinates not changed by the
$G_a$-action). Thus $X$ has infinite-dimensional automorphism group.
The same proof works for the complement $\P^{n-1}-\{q=0\}$ of
an isotropic quadratic form $q$, since there is a nontrivial homomorphism
from $G_a$ to the orthogonal group of $q$.
\qed

In particular, every quadratic form of dimension at least 2 over
an algebraically closed field is isotropic. For example, 
Lemma \ref{inf} says that the sphere $x_0^2+\cdots +x_n^2=1$,
viewed as an affine variety over the complex numbers, has
infinite-dimensional automorphism group for $n$ at least 2. In the simplest
case $n=2$, we can change variables over the complex numbers
to $x_1x_2+x_3^2=1$, and then the proof of Lemma \ref{inf} gives
the well-known automorphisms
$$(x_1,x_2,x_3)\mapsto (x_1-f(x_2)^2x_2-2f(x_2)x_3,x_2,f(x_2)x_2+x_3)$$
for any polynomial $f$ in one variable.

For affine quadrics whose homogeneous part is anisotropic, it is not
known whether the automorphism group of $X$ as an affine variety
is just the orthogonal group;
Theorem \ref{witt} proves this in many cases. By contrast, we now show
that the group of birational automorphisms of a smooth quadric $Q$
of dimension at least 2 is always infinite-dimensional. (We define this
to mean
that there are unirational families of arbitrarily large dimension
over $k$ of birational automorphisms of $Q$, as in the above
definition for biregular automorphisms.)

\begin{lemma}
Let $Q$ be a smooth projective quadric of dimension at least 2
over a field $k$. Then the group of birational automorphisms of 
$Q$ over $k$ is infinite-dimensional.
\end{lemma}

{\bf Proof. }
Consider a general linear projection from the quadric
$Q^n\subset \P^{n+1}$ onto a linear space $\P^{n-r}$,
for any $1\leq r\leq n-1$. This gives a rational map $f$ from
$Q^n$ to $\P^{n-r}$ of which the general fiber is an
$r$-dimensional quadric.
This general fiber comes from a quadratic form $\rho$
of dimension $r+2$ over the function field
$l=k(x_1,\ldots,x_{n-r})$ of $\P^{n-r}$. The group
of birational automorphisms of $Q$ contains the group of
birational automorphisms of $Q$ over $\P^{n-r}$, which in turn
contains the orthogonal group $(PGO(\rho))(l)$,
using that $r\geq 1$. (The definition of $PGO$ is recalled in
section \ref{finite}.) The identity component of the
algebraic group $PGO(\rho)$
is reductive, and so it is unirational over $l$
\cite[Theorem 18.2]{Borel};
also, it has positive dimension over $l=k(x_1,\ldots,
x_{n-r})$. Thus, the group $(PGO(\rho))(l)$ gives an infinite-dimensional
group of birational automorphisms of $Q$ over $k$.
\qed

\section{Finite-dimensional automorphism groups}
\label{finite}

We now present some useful general results on isomorphisms
between open varieties. The proofs are easy consequences
of Abhyankar's lemma in birational geometry.
As an application, for any
affine quadric whose homogeneous part of degree 2 is anisotropic
with first Witt index equal to 1, the automorphism group
is essentially the orthogonal group. The application uses
Karpenko's theorem that a projective quadric with first Witt
index equal to 1 is not ruled.

Abhyankar's lemma is the following statement
\cite[Proposition 4]{Abhyankar};
see Koll\'ar \cite[Theorem VI.1.2]{Kollarbook}
for a recent exposition.

\begin{lemma}
\label{abh}
Let $\pi:Y\arrow X$ be a proper birational morphism of irreducible
schemes, with $Y$ normal and $X$ regular. Then every exceptional
divisor of $\pi$ is ruled over its image. That is, if $E$
is an irreducible divisor in $Y$ which does not map birationally
to its image $E'$ in $X$, then $E$ is birational over $E'$
to a scheme $W\times_{E'}\P^1_{E'}$.
\end{lemma}

We deduce the following information about isomorphisms between
open varieties, of which part (2) is especially powerful.
A variety over a field $k$ is {\it ruled }if
it is birational to $Y\times \P^1$ for some variety $Y$ over $k$.

\begin{theorem}
\label{ruled}
(1) Let $(X_1,D_1)$ and $(X_2,D_2)$ be pairs with $X_i$ a projective
variety over a field $k$ and $D_i$ an irreducible divisor, for $i=1$ and 2.
We assume that each $X_i$ is regular in a neighborhood of $D_i$.
Suppose that $D_1$ is not ruled over $k$. Then any isomorphism $f$
from $X_1-D_1$ to $X_2-D_2$, viewed as a birational map
from $X_1$ to $X_2$, is an isomorphism in codimension 1.
That is, there are subsets $S_i$
of codimension at least 2 in $X_i$ such that $f$ is an isomorphism
from $X_1-S_1$ to $X_2-S_2$.

(2) Suppose in addition that the divisors $D_1$ and $D_2$ are ample.
Then any isomorphism
from $X_1-D_1$ to $X_2-D_2$ extends to an isomorphism from $X_1$
to $X_2$.
\end{theorem}

Under the stronger assumptions that
$D_1$ is not uniruled and the base field is the complex numbers,
essentially the same conclusion was proved by Jelonek
\cite[Theorem 3.7]{Jelonek}.
For our application, it is
crucial that we only need to assume that $D_1$ is not ruled:
every quadric over a field is uniruled, but many of them are not ruled.

{\bf Proof. }
(1) Let $Y$ be the normalization of the closure
of the graph of $f$ in $X_1\times X_2$; then we have proper birational
morphisms $f_i:Y\arrow X_i$ with $f=f_2f_1^{-1}$.
The proper transform of $D_1$ in $Y$
is a divisor birational to $D_1$. Since $D_1$ is not
ruled, Lemma \ref{abh} shows that $D_1$ maps birationally
to its image in $X_2$. Since $D_2$ is irreducible, this means that
$D_1$ maps birationally to $D_2$. Thus, $f$ induces an isomorphism
between some open subsets $X_i-S_i$ where $S_i$ has codimension
at least 1 in $D_i$, thus codimension at least 2 in $X_i$.

(2) The birational map $f$ from $X_1$ to $X_2$ is an isomorphism
in codimension 1, and it maps the divisor $D_1$ to the divisor $D_2$.
So $f$ induces an isomorphism between the rings
$\oplus_{j\geq 0}H^0(X_i, jD_i)$. (The point is that the group
of sections of a line bundle on a normal variety remains unchanged
upon removal of a subset of codimension at least 2 from the variety
\cite[Lemma 2.32]{Iitaka}.
We are assuming that $X_i$ is regular,
hence normal, in a neighborhood of $D_i$.)
Since $D_i$ is ample
for $i=1$ and 2, the projective variety Proj of that graded ring
is isomorphic to $X_1$ or $X_2$, respectively. So $f$ induces
an isomorphism from $X_1$ to $X_2$.
\qed

The same proof works if we 
replace the assumption that the divisors $D_i$ are ample
by the assumption that the varieties $X_i$ are Fano, that is,
that their anticanonical line bundles are ample. Both conditions
will be satisfied in our application, Theorem \ref{witt}.

We draw the following conclusion about automorphism groups
of affine quadrics. Define the first Witt index $i_1(q)$
of an anisotropic quadratic form $q$ of dimension at least 2
over a field to be the
Witt index (the maximum dimension of an isotropic subspace)
of the quadratic form $q$ over the function field of the projective
quadric $\{q=0\}$ \cite{Knebusch1}. We always have $i_1(q)\geq 1$.
If $V$ denotes the vector space on which $q$ is defined,
we write $GO(V)$ for the subgroup of linear automorphisms that preserve
$q$ up to scalars and $PGO(V)$ for $GO(V)/k^*$. More generally,
for a linear subspace $W\subset V$, $PGO(W\subset V)$ denotes
the subgroup of $PGO(V)$ that maps $W$ into itself.

\begin{theorem}
\label{witt}
Let $f$ be a polynomial of degree at most 2 on a vector space $W$
over a field $k$,
and let $q$ be the homogeneous part of degree 2 of $f$. Assume
that the quadratic form $q$ is anisotropic and that the first Witt index
$i_1(q)$ is equal to 1. If $k$ has characteristic 2, assume that
the projective quadric $\{q=0\}$ over $k$ is a regular scheme.
By introducing an extra variable to
homogenize $f$, we get
a quadratic form on a vector space $V$ which restricts to $q$
on a hyperplane $W\subset V$.
Then the automorphism group of the affine quadric
$\{f=0\}$ is equal to $PGO(W\subset V)$.
Also, the automorphism group of the complement $\P^{n-1}-\{q=0\}$,
where $n$ is the dimension of $q$,
is equal to $PGO(V)$.
\end{theorem}

{\bf Proof. }
For an anisotropic quadratic form $q$ with first
Witt index 1 over $k$, Karpenko proved
that the associated projective quadric $Q$ is not ruled over $k$
\cite[Theorem 6.4]{Karpenko}.
Karpenko's theorem was extended to smooth quadrics in characteristic 2
by Elman, Karpenko, and Merkurjev \cite{EKM} and then to all
quadrics in characteristic 2 \cite{Totaroess}.

Both the affine quadric $\{f=0\}$ and the complement
$\P^{n-1}-\{q=0\}$ are complements of $Q$ in projective
varieties $X$. Also, $X$ is regular in a neighborhood of $Q$ 
since the Cartier divisor $Q$ is regular. Finally,
the divisor $Q$ is ample on $X$.
By Theorem \ref{ruled},
every automorphism of the affine varieties $\{f=0\}$ and
$\P^{n-1}-\{q=0\}$ extends to
an automorphism of the compactification. It is then easy
to compute the automorphism group precisely.
\qed

To see what Theorem \ref{witt} means more concretely, let us discuss
which quadratic forms have first Witt index equal to 1.
Briefly, this is the ``typical''
behavior of quadratic forms of any dimension over a sufficiently complicated
field. For example,
the generic form $t_1x_1^2+\cdots+t_nx_n^2$ over $k=k_0(t_1,\ldots,t_n)$
has first Witt index 1 \cite[Example 5.7]{Knebusch1}.
See section \ref{questions} to see
the meaning of the first Witt index for low-dimensional forms.
Also, over any field,
every anisotropic quadratic form whose
dimension has the form $2^a+1$ has first Witt index 1,
by Hoffmann \cite{Hoffmannwitt} (in characteristic not 2)
and Hoffmann-Laghribi \cite{HL} (in characteristic 2). (Hoffmann's
theorem in characteristic not 2 was later reproved using
Rost's degree formula \cite{MerkurjevRost} and also by studying Steenrod
operations on Chow groups of quadrics \cite[Corollary 75.8]{EKM}.)
In particular, the form
$x_1^2+\cdots +x_r^2$ over the real numbers has first Witt index 1
whenever $r$ is of the form $2^a+1$; this case was known earlier by
Knebusch \cite[Proposition 7.9, Example 7.10]{Knebusch2}.
(The form $x_1^2+\cdots
+x_r^2$ over the real numbers 
with $r$ not of the form $2^a+1$ has first Witt index
greater than 1, as explained at the end of section \ref{questions}.)
We deduce the following conclusion, proved for $n=2$
in \cite[Theorem 6.2]{Totarocomp}.

\begin{corollary}
\label{autsphere}
The automorphism group of the sphere
\[ S^n_{\R}=\{x_0^2+\cdots + x_n^2=1\},\]
as an affine algebraic variety over the real numbers,
is equal to the orthogonal group $O(n+1)$ whenever $n$
is a power of 2.
\end{corollary}

It is not known whether the automorphism group of $S^n_{\R}$ is
only the orthogonal group when $n$ is not a power of 2.

Corollary \ref{autsphere} is vaguely reminiscent of Wood's theorem that
when $n$ is a power of 2, any morphism from the affine variety
$S^n_{\R}$ to a lower-dimensional sphere is constant \cite{Wood}.
It is interesting to note that Wood's theorem may not be optimal
in high dimensions. It is optimal in low dimensions, since the Hopf maps
are nonconstant polynomial maps $S^3\arrow S^2$, $S^7\arrow S^4$,
and $S^{15}\arrow S^8$. The first open case is: is there a nonconstant
polynomial map $S^{48}\arrow S^{47}$? At least there is no such map given
by homogeneous quadratic polynomials \cite{Yiu}.

\section{Which quadrics are ruled?}
\label{questions}

The heart of the proof of Theorem \ref{witt} is the fact that,
by Karpenko's theorem,
an anisotropic quadratic form with first Witt index 1
is not ruled (that is, the corresponding projective quadric is not ruled).
We conjecture the converse (Conjecture \ref{quadruled}). In this section
we give some evidence, including proving the conjecture for
quadratic forms of dimension at most 9 (Lemma \ref{nine}).

\begin{conjecture}
\label{quadruled}
Let $q$ be an anisotropic quadratic form over a field $k$.
If the first Witt index of $q$ is greater than 1, then $q$
is ruled over $k$.
\end{conjecture}

For the rest of this section, we will assume that the field $k$
is not of characteristic 2. The conjecture in characteristic 2 is
discussed in \cite{Totaroess}. Note that the problem of ruledness
is only interesting for anisotropic quadrics, since a smooth isotropic
quadric of dimension at least 1 over a field $k$
is rational over $k$ and hence ruled over $k$.

Some evidence for Conjecture \ref{quadruled} is that
a ``stabilized'' version is true. Namely, if $Q$ is an anisotropic quadric
with $i_1(Q)$ greater than 1, let $Q'$ be any subquadric of codimension
$r:=i_1(Q)-1$. Then the definition of $i_1(Q)$ implies that $Q'$
becomes isotropic over the function field $k(Q)$; that is, there
is a rational map from $Q$ to $Q'$ over $k$. Since there is also a
rational map from $Q'$ to $Q$ (the inclusion), a standard
result in the theory of quadratic forms gives that $Q$
is stably birational to $Q'$ \cite[Theorem X.4.25]{Lam}.
That is,
$Q\times \P^a$ is birational to $Q'\times \P^{r+a}$ for some $a$.
So it seems plausible (an analogue of the ``quadratic Zariski
problem'' in Ohm \cite{Ohm}) that $Q$ should be birational to
$Q'\times \P^r$ and hence that $Q$ should be ruled.

Conjecture \ref{quadruled} is true for special Pfister neighbors,
by Knebusch. We recall the definitions. We use Lam's book
\cite[Chapter X] {Lam}
as a reference on the theory of quadratic forms, although
our notation is slightly different. Two quadratic forms over a field $k$
are similar if one is isomorphic to a nonzero scalar multiple of the other.
An $n$-fold Pfister form $\langle \langle a_1,\ldots,a_n\rangle \rangle$
is the tensor product of the  2-dimensional
quadratic forms $\langle 1,-a_i\rangle$ for $a_1,\ldots,a_n$ in $k^*$.
A Pfister neighbor
is a form $q$ similar to a subform of an
$n$-fold Pfister form with
$\text{dim }q>2^{n-1}$. A special Pfister neighbor is 
a form $q$ similar to $\alpha\perp c\alpha'$ for some $(n-1)$-fold
Pfister form $\alpha$ and some nonzero subform $\alpha'$ of $\alpha$;
clearly $q$ is a neighbor of the $n$-fold Pfister form $\alpha\perp  c\alpha$.
One can check that a neighbor of an $n$-fold Pfister form is special
if and only if it contains a scalar multiple of some
$(n-1)$-fold Pfister form.
Pfister neighbors are among the simplest quadratic forms: in particular,
an anisotropic
Pfister neighbor of dimension $2^{n-1}+a$, where $1\leq a\leq 2^{n-1}$,
has first Witt index $a$, which is the largest possible for forms of that
dimension, by Hoffmann \cite{Hoffmannwitt}.

Thus, Conjecture \ref{quadruled} 
predicts that a Pfister neighbor whose dimension
is not of the form $2^{n-1}+1$ should be ruled. This is not known
in general, but it is true for special Pfister neighbors by
Knebusch; more precisely, he showed that the function field of
$\alpha\perp  c\alpha'$ is a purely transcendental extension of the
function field of the
``base form'' $\alpha\perp  \langle c \rangle$
\cite[pp.~73--74]{Knebusch1}.
Ahmad-Ohm observed that
Knebusch's argument produces other examples of ruled quadrics,
as follows \cite[after 1.3]{AO}.
In particular, every form divisible by a binary form is ruled.

\begin{lemma}
\label{pfister}
Let $P$ be a Pfister form over a field $k$, $P_1$ a nonzero subform of $P$,
and $b_1,\ldots, b_r$ any elements of $k^*$. Then the projective quadric
associated to $b_1P\perp \cdots\perp b_{r-1}P\perp b_rP_1$ is birational
to the product of a projective space with the quadric
$b_1P\perp \cdots\perp b_{r-1}P\perp \langle b_r \rangle$;
in particular, it is ruled if $P_1$ has dimension at least 2.
\end{lemma}

{\bf Proof. }
A Pfister form $P$ is strongly multiplicative: there is
a ``multiplication''
$xy$ on the vector space of $P$ which is a rational function of $x$ and
linear in $y$ such that $P(xy)=P(x)P(y)$ \cite[Theorem X.2.11]{Lam}.
Given a general point
on the first quadric $b_1P(x_1)+\cdots+b_{r-1}P(x_{r-1})+b_rP(x_r)=0$
(where $x_r\in P_1$), we map it to a point on the second quadric
by noting that $b_1P(x_rx_1)+\cdots +b_{r-1}P(x_rx_{r-1})+b_rP(x_r)^2=0$.
The general fibers of this rational map are linear spaces. \qed

\begin{lemma}
\label{nine}
Let $q$ be an anisotropic quadratic form
over a field $k$ of characteristic not 2. If the first Witt index of $q$
is greater than 1 and $q$ has dimension at most 9, then $q$ is ruled.
\end{lemma}

{\bf Proof.} Every form $q$ of dimension 3 or 5 has $i_1(q)=1$, and so
there is nothing to check. A form $q$ of dimension 4 with $i_1(q)>1$
has $i_1(q)=2$ and is similar
to a Pfister form \cite[Theorem X.4.14]{Lam}.
So the corresponding quadric surface is ruled.
A form $q$ of dimension 6 with $i_1(q)>1$ has $i_1(q)=2$
and is divisible by a binary form, by Knebusch
\cite[Theorem 10.3]{Knebusch2}.
Therefore $q$ is ruled.
A form $q$ of dimension 7 with $i_1(q)>1$ has $i_1(q)=3$ 
\cite[Theorem 4.1]{Hoffmannsplit}
and is therefore similar to the pure
subform $p'=p-\langle 1\rangle$ of a Pfister form $p$
\cite[Theorem 5.8]{Knebusch1}.
Clearly
such a Pfister neighbor is special, and so $q$ is ruled.
Next, if $q$ is a form of dimension 8
with $i_1(q)>1$, then either $i_1(q)=4$ and $q$ is similar
to a 3-fold Pfister form, or $i_1(q)=2$ and $q$ is divisible by
a binary form,
by Hoffmann \cite[Theorem 4.1]{Hoffmannsplit}.
It follows that $q$ is ruled.
Finally, every form $q$
of dimension 9 has $i_1(q)=1$ and so there is nothing to prove. \qed

For 10-dimensional forms, Conjecture \ref{quadruled} remains open.
By Izhboldin, a 10-dimensional form $q$ with $i_1(q)>1$ is either
divisible by a binary form or is a Pfister neighbor
\cite[proof of Conjecture 0.10]{Izhboldin}.
In the first case,
$q$ is ruled, but the second case is harder.
There are
10-dimensional Pfister neighbors which are not special, by Ahmad-Ohm
\cite[2.8]{AO},
and it is not clear how to prove that such forms
are ruled.

Another example: for any $n$, the form
$x_0^2+\cdots +x_n^2$ is a special Pfister neighbor over any
field. That is clear by viewing this form
as a subform of a Pfister form
$\langle\langle  -1,\cdots, -1\rangle\rangle$.
So, for example over the real numbers, the
quadric $\{x_0^2+\cdots+x_n^2=0\}$ is ruled if and only if $n$ is not a
power of 2, in agreement with Conjecture \ref{quadruled}.

{\bf Acknowledgement.}
Thanks to Jean-Louis Colliot-Th\'el\`ene and Detlev Hoffmann for
their comments.


\end{document}